\documentclass[10pt]{IEEEtran}
\usepackage{graphicx}
\usepackage{color}
\usepackage{dsfont}
\usepackage{amsmath}
\usepackage{amssymb}
\usepackage{cite}
\usepackage{float}
\usepackage{stfloats}
\usepackage{bbm}
\usepackage{url}
\usepackage[geometry]{ifsym}

\newcommand\sups[1]{^{\hbox{\scriptsize{#1}}}}
\newcommand\supt[1]{^{\hbox{\tiny{#1}}}}
\newcommand\subs[1]{_{\hbox{\scriptsize{#1}}}}
\newcommand\subt[1]{_{\hbox{\tiny{#1}}}}

\newcommand{\vb}[1]{\mathbf{#1}}
\newcommand{\mb}[1]{\mathbb{#1}}
\newcommand{\numeq}[2]{\begin{equation} #2 \label{#1} \end{equation}}

\newcommand{\mc}[1]{\mathcal{#1}}

\newcommand{\citeasnoun}[1]{Ref.~\citenum{#1}}

\newcommand{\tet}{\text{\TriangleUp}}
\newcommand{\tetbar}{\overline{\text{\TriangleUp}}}
\newcommand{\Tet}{\mc T}

\newcommand{\vbXi}{\boldsymbol{\xi}}
\newcommand{\vbEta}{\boldsymbol{\eta}}
\newcommand{\citeasnouns}[2]{Refs.~\citenum{#1},~\citenum{#2}}

\begin{document}

\title {Taylor-Duffy Method for Singular Tetrahedron-Product Integrals:
        Efficient Evaluation of Galerkin Integrals for VIE Solvers} 
\author{M.~T.~Homer~Reid%
        \thanks{M. T. Homer Reid is with the Department of Mathematics,
                Massachusetts Institute of Technology.}
       }
\maketitle

\begin{abstract}

I present an accurate and efficient technique for
numerical evaluation of singular 6-dimensional
integrals over tetrahedon-product domains,
with applications to calculation of Galerkin matrix
elements for discretized volume-integral-equation (VIE) 
solvers using Schaubert-Wilton-Glisson (SWG) and other 
tetrahedral basis functions. My method
extends the generalized Taylor-Duffy strategy---used to
handle the singular \textit{triangle}-product integrals
arising in discretized surface-integral-equation (SIE)
formulations---to the tetrahedron-product case; it 
effects an exact transformation of a singular 6-dimensional
integral to an nonsingular lower-dimensional integral
that may be evaluated by simple numerical cubature
The method is highly general and may---with the
aid of automatic code generation facilitated by
computer-algebra systems---be applied to a wide
variety of singular integrals arising in various
VIE formulations with various types of tetrahedral
basis function, of which I present several examples.
To demonstrate the accuracy and efficiency of my method,
I apply it to the calculation of matrix elements
for the volume electric-field integral equation (VEFIE)
discretized with SWG basis functions, where the
method yields 12-digit or higher accuracy with
low computational cost---an improvement of
many orders of magnitude compared to existing
techniques.
\end{abstract}

%------------------------------------------------------------
%- Document body
%------------------------------------------------------------

\section{Introduction}

In this paper I present an efficient technique for evaluating
singular 6-dimensional integrals over tetrahedron-product
domains, such as those commonly encountered in discretized
volume-integral equation (VIE) 
formulations~\cite{Harrington93, Chew2009, Volakis2012}
with tetrahedral basis 
functions~\cite{Bleszynski2008, Markkanen2014, Schaubert1984}.
My method extends the generalized Taylor-Duffy method for
singular triangle-product integrals~\cite{Duffy1982, Taylor2003, Reid2015A}
to the tetrahedron-product case, exactly transforming
singular 6-dimensional integrals to nonsingular lower-dimensional
integrals amenable to simple numerical cubature. I formulate the
basic algorithm, show how it may be applied to several
distinct VIE formulations,
and present computational results demonstrating its
accuracy and efficiency; for the specific case
of singular VEFIE integrals~\cite{Chew2009} with SWG basis
functions~\cite{Schaubert1984} (defined below)
I obtain 12 or more digits of accuracy with modest
computational cost, an improvement of many orders
of magnitude compared to a recently-proposed
alternative approach~\cite{Bleszynski2013}.

Discretized VIE methods~\cite{Harrington93, Chew2009, Volakis2012}
using tetrahedral basis functions~\cite{Schaubert1984,Markkanen2014}
are useful for attacking many problems in science and
engineering, including electromagnetic 
scattering~\cite{Hasanovic2007, Zhang2015},
acoustic wave propagation~\cite{Bleszynski2008},
inductance extraction~\cite{Jackman2016},
and fluctuation-induced phenomena~\cite{Reid2016B}.
Although there exists a considerable variety of
VIE formulations~\cite{Botha2006,Bleszynski2008,Markkanen2012, Polimeridis2014}
and multiple choices of tetrahedral basis
functions---including piecewise-constant~\cite{Bleszynski2008},
piecewise-linear~\cite{Markkanen2014}, and
SWG functions~\cite{Schaubert1984}---a computational
challenge common to \textit{all} Galerkin VIE formulations
is the need for accurate and efficient numerical 
evaluation of 6-dimensional tetrahedron-product integrals,
typically of the general form
%====================================================================%
\numeq{OriginalIntegral}
{
 \mathcal{I}
=\int_{\Tet} d\vb x \, \int_{\Tet^\prime} d\vb x^\prime \, 
  P(\vb x, \vb x^\prime) K(|\vb x-\vb x^\prime|)
}
%====================================================================%
where $\Tet,\Tet^\prime$ are tetrahedra,
$P$ is a polynomial in the cartesian components of
$\vb x, \vb x^\prime$, and $K(r)$ is a scalar kernel function.
In the commonly-encountered case in which $\Tet,\Tet^\prime$
have one or more common vertices and $K(r)$ is singular
at $r=0$, the integral (\ref{OriginalIntegral}) cannot be
evaluated by simple numerical cubature~\cite{Cools2003};
instead, more sophisticated integration strategies
are required, whose accuracy and efficiency play a large
part in determining those of the overall VIE 
solver~\cite{Botha2006}.

To date, several methods for evaluating singular 
tetrahedron-product integrals have been discussed;
many such methods, including that proposed here,
extend techniques originally developed for singular
\textit{triangle}-product integrals [like (\ref{OriginalIntegral})
but with $\Tet,\Tet^\prime$ replaced by two-dimensional triangular 
domains], a problem that has been studied for decades
due to its importance for surface-integral-equation
(SIE) solvers~\cite{RWG1982, Medgyesi1994, Chew2009}.
Strategies proposed for singular tetrahedon-product integrals
include singularity subtraction~\cite{Hu2008, Jarvenpaa2003},
separation of inner and outer 3D integrals with the former
(latter) evaluated analytically (numerically)~\cite{Wilton1984},
and the use of Stokes' theorem~\cite{DFN1985} to recast volume 
integrals as surface integrals~\cite{Bleszynski2013}.

The method of Bleszynski et al.~\cite{Bleszynski2013}
is particularly attractive in that it effects an exact
transformation of the singular 6-dimensional integral
(\ref{OriginalIntegral}) to a sum of non-singular 
4-dimensional integrals amenable
to straightforward low-order numerical cubature.
However, the method was presented in~\citeasnoun{Bleszynski2013}
only for one particular VIE formulation with particular basis 
functions \{namely, the volume electric-field integral equation 
(VEFIE)~\cite{Markkanen2012, Polimeridis2014}
with SWG functions~\cite{Schaubert1984}\}, and it is unclear
if or how the method could be used to evaluate
the integrals arising in formulations.
Moreove, the number of 4-dimensional integrals
that must be evaluated is large (as many as $16$), and
these integrals---though nonsingular---converge
relatively slowly in numerical quadrature schemes,
a point that was noted already in \citeasnoun{Bleszynski2013} 
and which I corroborate and discuss in further 
detail below (see Section \ref{ComputationalResultsSection}
and Figure \ref{ConvergenceFigure}).

In this paper I propose a strategy that, like that
of~\citeasnoun{Bleszynski2013}, exactly transforms
(\ref{OriginalIntegral}) into a sum of non-singular
lower-dimensional integrals, which are evaluated by
low-order numerical cubature in a practical solver.
However, in contrast to the Stokes'-theorem underpinning
of~\citeasnoun{Bleszynski2013}, my method is based
on Duffy's singularity-cancellation
technique~\cite{Duffy1982},
which was applied to triangle-product integrals originally by
Taylor~\cite{Taylor2003} and later in more generality
by~\citeasnoun{Reid2015A}; here (Section II) I extend
these ideas to the tetrahedron-product case, culminating
in a nonsingular reduced-dimensional integral
[equation (\ref{FinalIntegral})] that is
exactly equivalent to (\ref{OriginalIntegral}).
This reduction scheme offers several advantages (Section III):
\textbf{(a)} I formulate the algorithm in full
generality for \textit{any} integral of the form
(\ref{OriginalIntegral}), offering immediate
application to many VIE formulations and choices
of basis function; some explicit examples are 
given in Section \ref{VIEFormulationsSection}.
\textbf{(b)} The dimension of the reduced integral
produced by this method is $D=6-N\subt{CV}$ in general
and $D=5-N\subt{CV}$ for power-law kernels $[K(r)\sim r^p]$,
where $N\subt{CV}=\{1,2,3,4\}$ is the number of vertices common
to $\Tet,\Tet^\prime$; in the common-tetrahedron $(N\subt{CV}=4)$
and common-face ($N\subt{CV}=3$) cases this results in 
final integrals of dimension $D=2$ or $D=3$
($D=1$ or $2$ for power-law kernels), a further
reduction of dimension than is achieved by the 
method of~\citeasnoun{Bleszynski2013}.
\textbf{(c)} The number of $D$-dimensional integrals
into which (\ref{OriginalIntegral}) is transformed
is at most 18 and as few as 9 in some cases; moreover,
all integrals extend over the same region of 
integration---namely, the unit $D$-dimensional 
hypercube---and may thus be combined into a single
integral, affording significant efficiency through 
reuse of computation.
Thus my method is not only quite general but also
highly efficient, as I demonstrate with
illustrative computational results in 
Section \ref{ComputationalResultsSection}.
Questions for future work are discussed in 
Section \ref{ConclusionsSection}, and 
technical details are relegated to Appendices.

The algorithm of this paper is implemented in 
{\sc buff-em} a free, open-source
software implementation of the VEFIE with
SWG basis functions~\cite{buff-em}.

\section{Extension of Taylor-Duffy Method to Tetrahedron-Product
Integrals}
\label{DerivationSection}

The method of Duffy transforms~\cite{Duffy1982} was
applied to the desingularization and dimensional reduction
of singular triangle-product integrals by
Taylor~\cite{Taylor2003} and later in more generality
by~\citeasnoun{Reid2015A}.
In this section I show that the same basic ideas may be
used to desingularize and reduce the dimension of
tetrahedron-product integrals; the result is
equation (\ref{FinalIntegral}), a
nonsingular reduced-dimensional integral
that is exactly equivalent to the singular
six-dimensional integral (\ref{OriginalIntegral}).

The logical flow of the transformation procedure is
identical to that of~\citeasnouns{Taylor2003}{Reid2015A}
and proceeds as follows.
\textbf{(a)} Subdivide the tetrahedron-product
domain in (\ref{OriginalIntegral}) into $D$
subdomains and change integration variables
to ensure that each subdomain is a product of tetrahedra
with one vertex at the origin, facilitating 
Duffy transformation
[Section (\ref{DerivationStep1})]. (For the
tetrahedron-product case we have $D=18$,
in contrast to $D=6$ for the
triangle-product case~\cite{Taylor2003,Reid2015A}.)
\textbf{(b)} Within each subdomain, analytically
evaluate the integrals over all variables of which the
kernel $K(r)$ is independent [Section (\ref{DerivationStep2})].
(As in the triangle-product case, there are
$N\subt{CV}-1$ such variables, where $N\subt{CV}$
is the number of vertices common to $\Tet,\Tet^\prime$.)
\textbf{(c)} Within each subdomain, perform
a Duffy transformation, analytically evaluate the
integral over the untransformed variable,
then combine the remaining integrals for each
subdomain into a single integral over the
$(6-N\subt{CV})$-dimensional unit hypercube
[Section (\ref{DerivationStep3})] to yield
the final master formula (\ref{FinalIntegral}).

As in the triangle-product case, the reduction
procedure is straightforward but tedious and
error-prone if carried out by hand, in practice
requiring the use of automatic code generation
facilitated by computer algebra systems
(Section \ref{ComputerAlgebraSection}).

Although the TD reduction method may be used
for all tetrahedron pairs with $N\subt{CV}=1$
or more common vertices,
here I formulate it only for the case of $N\subt{CV}\ge 2,$
as the significant cost of implementing the
TD method seem not to be justified by the modest 
reduction in computational cost it affords in the
common-vertex case.

%====================================================================%
%====================================================================%
%====================================================================%
\subsection{Decomposition into tetrahedron-product subdomains}
\label{DerivationStep1}

The goal of this step is to decompose the tetrahedron-product
domain in (\ref{OriginalIntegral}) as the union of $D$
tetrahedron-product subdomains,
%====================================================================%
$$ \Tet \times \Tet^\prime
   = \displaystyle{\bigcup_{d=1}^{D}} \tet_d\times\tetbar_d,
$$
%====================================================================%
with the property that $\tet_d$ and $\tetbar_d$ each
have one vertex at the origin of coordinates, as required
to allow Duffy transformation.

In (\ref{OriginalIntegral}) first make the change of variables
$(\vb x, \vb x^\prime) \to (\vbXi, \vbEta)$
where $\vbXi,\vbEta$ run over a standard tetrahedron $\tet_0$: 
% with vertices $\{(0,0,0), (1,0,0), (1,1,0), (1,1,1)\}:$
%%%%%%%%%%%%%%%%%%%%%%%%%%%%%%%%%%%%%%%%%%%%%%%%%%%%%%%%%%%%%%%%%%%%%%
\numeq{OriginalIntegral2}
{ \mathcal{I}= J\int_{\tet_0} \, d\vbXi \, \int_{\tet_0} d\vbEta \,
  P\Big( \vb x(\vbXi), \vb x^\prime(\vbEta)\Big)
  K\Big( r(\vbXi, \vbEta)\Big)
}
%%%%%%%%%%%%%%%%%%%%%%%%%%%%%%%%%%%%%%%%%%%%%%%%%%%%%%%%%%%%%%%%%%%%%%
where
%%%%%%%%%%%%%%%%%%%%%%%%%%%%%%%%%%%%%%%%%%%%%%%%%%%%%%%%%%%%%%%%%%%%%%
$$\int_{\tet_0} d\vbXi =
  \int_0^1 d\xi_1 \, 
  \int_0^{\xi_1} d\xi_2 \, 
  \int_0^{\xi_2} d\xi_3.
$$
and
%%%%%%%%%%%%%%%%%%%%%%%%%%%%%%%%%%%%%%%%%%%%%%%%%%%%%%%%%%%%%%%%%%%%%%
\begin{align*}
 \vb x(\vbXi) &= V_1 + \xi_1 \vb L_1
              + \xi_2 \vb L_2
              + \xi_3 \vb L_3
\\
 \vb x^\prime(\vbEta) &= V_1 + \eta_1 \vb L_1^\prime
                     + \eta_2 \vb L_2^\prime
                     + \eta_3 \vb L_3^\prime
\\
 r(\vbXi, \vbEta)&=\big|\vb x(\vbXi) - \vb x^\prime(\vbEta)\big|
\end{align*}
with
%%%%%%%%%%%%%%%%%%%%%%%%%%%%%%%%%%%%%%%%%%%%%%%%%%%%%%%%%%%%%%%%%%%%%%
\begin{align*}
 \vb L_1&=(\vb V_2-\vb V_1), \qquad 
 \vb L_2=(\vb V_3-\vb V_2),  \qquad
 \vb L_3=(\vb V_4-\vb V_3)
\\
 \vb L_1^\prime&=(\vb V_2^\prime-\vb V_1),       \qquad 
 \vb L_2^\prime=(\vb V_3^\prime-\vb V_2^\prime), \qquad
 \vb L_3^\prime=(\vb V_4^\prime-\vb V_3^\prime)
\end{align*}
%%%%%%%%%%%%%%%%%%%%%%%%%%%%%%%%%%%%%%%%%%%%%%%%%%%%%%%%%%%%%%%%%%%%%%
and $J=36 \mc V \mc V^\prime$.
Here $\{\vb V_i, \vb V^\prime_i\}$ and $\mc V, \mc V^\prime$
are the vertices and volumes of $\tet, \tet^\prime$.
We have assumed that $\mb T, \mb T^\prime$ have at least
one common vertex, labeled $\vb V_1$, and the remaining
vertices should be ordered such that common vertices
have lower indices than non-common vertices; thus
for the common-edge case we have $\vb V_2^\prime=\vb V_2$,
for the common-triangle case we have additionally 
$\vb V_3^\prime=\vb V_3$,
and for the common-tetrahedron case we have additionally
$\vb V_4^\prime=\vb V_4$.
Note that $r(\vbXi, \vbEta)$ is the square root of a homogeneous
second-degree polynomial in the components of $\vbXi, \vbEta$:
%%%%%%%%%%%%%%%%%%%%%%%%%%%%%%%%%%%%%%%%%%%%%%%%%%%%%%%%%%%%%%%%%%%%%%
\numeq{rDef}
{
 r(\vbXi,\vbEta)=
   \sqrt{ R_{ij}^{\xi\xi}\xi_i \xi_j
         +R_{ij}^{\xi\eta}\xi_i \eta_j
         +R_{ij}^{\eta\eta}\eta_i \eta_j
        }
}
%%%%%%%%%%%%%%%%%%%%%%%%%%%%%%%%%%%%%%%%%%%%%%%%%%%%%%%%%%%%%%%%%%%%%%
where the $\{R_{ij}\}$ coefficients are functions of the geometric
parameters; for example, $R_{12}^{\xi\xi}=2\vb L_1 \cdot \vb L_2.$
The homogeneity of the polynomial under the radical,
which follows from the fact that $\vb T, \vb T^\prime$ have
one or more common vertices, is what allows analytical
evaluation of the $w$ integral below 
(Section \ref{DerivationStep3}).

Following~\citeasnoun{Taylor2003}, now introduce the relative
coordinates $\vb u \equiv \vbEta - \vbXi$, change variables
from $\{\vbXi, \vbEta\}$ to $\{\vbXi, \vb u\}$, and decompose
the domain of integration into $D$ subdomains with
the property that, within each subdomain, both $\vb u$
and $\vbXi$ run over tetrahedra with one vertex at the origin:
%%%%%%%%%%%%%%%%%%%%%%%%%%%%%%%%%%%%%%%%%%%%%%%%%%%%%%%%%%%%%%%%%%%%%%
\numeq{OriginalIntegral3}
{ \mc I=\sum_{d=1}^{D}
   \int_{\tet_d^{\vb u}} \, d\vb u
   \int_{\tet_d^{\vbXi}} \, d\vbXi
    P\big( \vb x(\vbXi), \vb x^\prime(\vbXi+\vb u) \big)
    K\big( r(\vbXi, \vbXi+\vb u) \Big).
}
%%%%%%%%%%%%%%%%%%%%%%%%%%%%%%%%%%%%%%%%%%%%%%%%%%%%%%%%%%%%%%%%%%%%%%
The corresponding step in the triangle-product 
case~\cite{Taylor2003,Reid2015A}
similarly writes the original triangle-product integral
as a sum of integrals over triangle-product subdomains,
with both triangles in each subdomain having one vertex at
the origin~\cite{Taylor2003}. However, whereas that case
involves $D=6$ triangle-product subdomains,
for the tetrahedron-product case one finds that the minimum
number of subdomains allowing (\ref{OriginalIntegral})
to be decomposed in the form (\ref{OriginalIntegral3})
is $D=18.$ (For the common-tetrahedron case this
number may be reduced to $D=9$ by identifying
pairs of identical subdomains; this is analogous
to the reduction from $D=6$ to $D=3$ available
for the $N\subt{CV}=3$ case of the Taylor-Duffy
approach to triangle-product integrals~\cite{Taylor2003}).

Explicit definitions of the 18 tetrahedral subdomains
$\tet_d^{\vb u}, \tet_d^{\vbXi}$ are given in
the Appendix (Tables \ref{TetTetBarTable1},
\ref{TetTetBarTable2}).

%====================================================================%
%====================================================================%
%====================================================================%
\subsection{Analytical evaluation of $\vbXi$ integrals}
\label{DerivationStep2}

If the original tetrahedra have 2 or more common vertices,
the distance function $r$ in (\ref{OriginalIntegral2})
is independent of one or more of the $\vbXi$ variables.
[For example, in the common-edge case $(N\subt{CV}=2)$ 
$r$ is independent of $\xi_1$, while in the
common-tetrahedron case $(N\subt{CV}=4)$
$r$ is independent of all $\vbXi$ variables.]
The kernel factor $K(r)$ in (\ref{OriginalIntegral3}) may then
be pulled out of the integrals over those variables,
leaving integrals over just the polynomial $P$; these
may be evaluated analytically to yield new polynomials
$\overline{P}$ depending on just the remaining variables:
%=================================================
\begin{subequations}
\begin{align}
\overline{P}\supt{4CV}_d(\vb u_d)
&\equiv 
 \int \, d\xi_3 
 \int \, d\xi_2 
 \int \, d\xi_1 
 P(\vbXi, \vb u_d + \vbXi)
\\
%--------------------------------------------------------------------%
\overline{P}\supt{3CV}_d(\vb u_d, \xi_3)
&\equiv 
 \int \, d\xi_2 
 \int \, d\xi_1
 P(\vbXi, \vb u_d + \vbXi)
\\
%--------------------------------------------------------------------%
\overline{P}\supt{2CV}_d(\vb u_d, \xi_2, \xi_3)
&\equiv 
 \int \, d\xi_1
 P(\vbXi, \vb u_d + \vbXi).
\end{align}
\label{PBarDef}
\end{subequations}
%=================================================
This reduces the dimension of the integral by
$N\subt{CV}-1$.

%=================================================
%=================================================
%=================================================
\subsection{Duffy Transformation and evaluation of $w$ integral}
\label{DerivationStep3}

For each of the $D$ subregions I now make
a Duffy transformation~\cite{Duffy1982}---that is,
for $1\le d\le D $ I introduce functions
%====================================================================%
\numeq{DuffyTransformation}
{ \vb u_d=\vb u_d(w,\vb y), \qquad \vbXi_d=\vbXi_d(w,\vb y)}
%====================================================================%
and make the change of variables
$(\vb u_d,\vbXi_d)\to (w,\vb y)$
in the $d$th subregion integral;
here the dimension of the $\vb y$ vector is
$Y\equiv 6-N\subt{CV}$.
The Duffy transformations for each subregion
are tabulated in the Appendix
(Tables \ref{DuffyTableNCV4}-\ref{DuffyTableNCV2}).

As in the triangle-product case~\cite{Taylor2003, Reid2015A},
the key property of this transformation is that,
when expressed as functions of the new variables,
each component of $\vb u$ and $\vbXi$ is
proportional to $w$. 
This yields a Jacobian factor for the $d$th 
subdomain of the form 
$J_d(w,\vb y)\equiv w^{Y} \mc J_d(\vb y)$
and---in view of the homogeneity of the
polynomial in (\ref{rDef})---allows the quantity
$w$ to be extracted
from the square root in equation (\ref{rDef}):
%=================================================
\numeq{XdDef}
{r(\vbXi, \vb u) = w X_d(\vb y)}
%=================================================
with $X_d(\vb y)$ \textit{nonvanishing} over
the region of integration.
Also, the $\overline{P}$ polynomials defined
by (\ref{PBarDef}) may be expanded as power
series in $w$, with $\vb y$-dependent coefficients:
%=================================================
\numeq{ScriptPDef}
{ \overline{P}_d(\vb u, \vbXi)
   \equiv \sum_n \mc{P}_{dn}(\vb y) w^n
}
%=================================================
Finally, because the domain of integration for the 
Duffy-transform variables $(w,\vb y)$ is the same for 
all $D$ subdomains (namely, $0\le w, y_i \le 1$),
the order of summation and integration
in (\ref{OriginalIntegral2}) may be reversed
to yield a single integral whose integrand is a 
sum of $D$ terms;
using (\ref{XdDef}) and (\ref{ScriptPDef}),
the final transformed version of the original 
integral (\ref{OriginalIntegral}) then reads
%=================================================
\numeq{FinalIntegral}
{ \mathcal{I} = \int_{\Box^Y}
   \sum_{d=1}^D \mc J_d(\vb y)
   \sum_{n} \mc P_{dn}(\vb y) \mc K_{n+Y}\Big(X_d(\vb y)\Big)
   \, d\vb y
}
%=================================================
where the $\{\mc K\}$ functions are the ``first integrals''
of $K$, defined by~\cite{Reid2015A}
%=================================================
\numeq{FirstIntegral}
{\mc K_p(X) \equiv \int_0^1 w^p K(wX) dw.}
%=================================================
If the kernel $K(r)$ has a singularity of degree $q$
at the origin [i.e. $K(r)\sim \frac{1}{r^q}$ as $r\to 0$]
then $\mc K_p$ exists and is nonsingular for $p\ge q$;
equation (\ref{FinalIntegral}) thus desingularizes
all integrals of the form (\ref{OriginalIntegral})
with singularities as strong as 
$\{\frac{1}{r^2}, \frac{1}{r^3}, \frac{1}{r^4}$ for
the common-\{tetrahedron, face, edge\} cases. If
the polynomial $P(\vb x,\vb x^\prime)$ vanishes 
at $\vb x=\vb x^\prime$, then the sum over $n$
in (\ref{FinalIntegral}) begins at $n=1$ or higher,
in which case kernels with even
stronger singularities are desingularized by 
(\ref{FinalIntegral}); an example is given in 
Section \ref{VIEFormulationsSection}.

As noted in~\citeasnouns{Taylor2003}{Reid2015A},
the first integral (\ref{FirstIntegral}) may be
evaluated in closed form for many kernels relevant to VIE
solvers, including the Helmholtz kernel
$\frac{e^{ikr}}{4\pi r}$ and its gradient;
explicit expressions for $\mc K_n$ for various 
kernels of interest may be found in~\cite{Reid2015A}.
For some kernels---in particular, power-law kernels
of the form $K(r)\sim r^p$ for integer $p$---the
dimension of the transformed integral (\ref{FinalIntegral})
may be \textit{further} reduced by evaluating one of the $\vb y$
integrals analytically, yielding a final reduced
integral of dimension $\{1,2,3\}$ for the 
common-\{tetrahedron, face, edge\} cases. This is
useful for applications to singularity-subtraction
methods~\cite{Jarvenpaa2003} or frequency-caching
schemes~\cite{Reid2015A} in which the contributions
of individual terms in the power-series expansion
of $K(r)$ around $r=0$ are computed analytically,
yielding integrals of the form (\ref{OriginalIntegral})
with $K(r)\sim r^p$ for various powers $p$.
The procedure for effecting this further reduction
of (\ref{FinalIntegral}) is identical to that 
presented in~\citeasnoun{Reid2015A} 
for the triangle-product case.

On the other hand, even for cases in which the
kernel is so complicated that even
the first integrals (\ref{FirstIntegral})
cannot be evaluated in closed form, equation
(\ref{FinalIntegral}) together with (\ref{FirstIntegral})
may still constitute a useful exact transformation of the
original integral (\ref{OriginalIntegral});
if the original integral has integrable singularities
that prevent direct application of numerical cubature,
equations (\ref{FinalIntegral}) and (\ref{FirstIntegral})
define a non-singular integral to which straightforward
numerical cubature may be applied directly.

%====================================================================%
%====================================================================%
%====================================================================%
\subsection{Automation by computer algebra system}
\label{ComputerAlgebraSection}

The reduction procedure outlined above,
though conceptually straightforward, in practice
requires large numbers of elementary calculus and 
algebra manipulations that are tedious and
error-prone if carried out by hand. Indeed,
to obtain the integrand of the reduced integral
(\ref{FinalIntegral}) for a given polynomial
$P$ and kernel $K$ we must---for each
of the $D$ subregions---\textbf{(a)} evaluate the 
integrals in (\ref{PBarDef}) to compute
the functions $\overline{P}_d(\vb u, \vbXi)$,
\textbf{(b)}
use the $d$-dependent Duffy transformation
(\ref{DuffyTransformation}) to
rewrite in terms of $(w,\vb y)$,
\textbf{(c)}
series-expand in $w$ to identify the 
coefficient functions $\mc P_{dn}(\vb y)$
in (\ref{ScriptPDef}), then 
\textbf{(d)}
pair each $\mc P$ with the appropriate
$\mc J$ and $\mc K$ factors and sum over
subdomains $d$ to construct
a function of $\vb y$ that may be
passed to a numerical cubature routine
as the integrand of (\ref{FinalIntegral})
The resulting integrand routinely consists of 
hundreds of terms, with the complexity
increasing with that of the polynomial $P$ in
(\ref{OriginalIntegral}); attempts 
to construct this function by hand are 
clearly hopeless.

As in the triangle-product case~\cite{Reid2015A},
the solution is to make avail to 
code generation by computer algebra systems
such as {\sc mathematica} or {\sc maxima},
which are ideally suited to carrying out
steps (\textbf{a}-\textbf{d}) above
automatically and emitting code
defining the integrand of (\ref{FinalIntegral}).
This approach was used to implement the
method of this paper in the 
{\sc buff-em} FVC solver~\cite{buff-em},
and the online documentation for that solver
includes sample {\sc mathematica} codes 
implementing the code-generation process
for particular cases of (\ref{OriginalIntegral}).

\section{Applications to Galerkin VIE Formulations with Tetrahedral Basis Functions}

The Taylor-Duffy method proposed in this paper is directly applicable to the
computation of matrix elements for a variety of Galerkin-discretized
VIE formulations using tetrahedral basis functions.
Here I first give the specific forms of the $P$ polynomial
and $K$ kernel in (\ref{OriginalIntegral}) for several popular
VIE formulations (Section \ref{VIEFormulationsSection}),
then present computational results for the particular 
case of the VEFIE formulation~\cite{Helsinki2014} with 
SWG basis functions~\cite{Schaubert1984}
(Section \ref{ComputationalResultsSection}).
As I show, the reduced integral (\ref{FinalIntegral}) produced by the
Taylor-Duffy transformation in this case may be evaluated to
high (12-digit or greater) accuracy with only $\sim 20$ quadrature
points per dimension, an improvement of many orders of magnitude
compared to another recently-proposed technique for singular
tetrahedron-product integrals~\cite{Bleszynski2013}.

\subsection{$P$ and $K$ functions for various VIE formulations}
\label{VIEFormulationsSection}

\textit{AIM acoustic-wave solver with tetrahedron-pulse functions.}
For the adaptive integral method (AIM) acoustic-wave VIE
formulation with piecewise-constant (pulse) tetrahedral basis
functions~\cite{Bleszynski2008}, elements of the stiffness matrix
take the form of (\ref{OriginalIntegral}) with the $P$ and $K$ 
functions given by
%====================================================================%
$$ P^\text{AIM}(\vb x, \vb x^\prime)=1,\qquad
   K^\text{AIM}(r) =\frac{e^{ikr}}{4\pi r}
$$
%====================================================================%
with $k$ the acoustic wavenumber in the background medium.

\textit{VEFIE with SWG functions.} For the volume electric-field
integral equation (VEFIE)~\cite{Helsinki2014} discretized with SWG basis
functions~\cite{Schaubert1984}, each element of the system matrix is a sum
of four tetrahedon-product integrals of the form
(\ref{OriginalIntegral}) with the $P$ and $K$ functions given by
%====================================================================%
\begin{subequations}
\begin{align}
 P\sups{EFIE}(\vb x, \vb x^\prime)
&=(\vb x - \vb Q) \cdot (\vb x^\prime - \vb Q^\prime) - \frac{9}{k^2},
\\
K\sups{EFIE}(r)
&=\frac{e^{ikr}}{4\pi r}
\end{align}
\label{PKEFIE}
\end{subequations}
%====================================================================%
Here $\vb Q, \vb Q^\prime$ are the source/sink vertices of the SWG
functions and $k$ is the vacuum photon wavenumber.

\textit{VMFIE with SWG functions.} For the volume magnetic-field
integral equation (VMFIE)~\cite{Helsinki2012} with SWG basis
functions, each element of the system matrix is a sum
of four tetrahedon-product integrals of the form
(\ref{OriginalIntegral}) with
%====================================================================%
\begin{subequations}
\begin{align*}
 P^\text{MFIE}(\vb x, \vb x^\prime)
&=(\vb x - \vb Q)\cdot
  \Big[ (\vb x-\vb x^\prime) \times (\vb x^\prime - \vb Q^\prime)\Big]
\\
K^\text{MFIE}(r)
&=(ikr-1)\frac{e^{ikr}}{4\pi r^3}.
\end{align*}
\label{PKMFIE}%
\end{subequations}
%====================================================================%
The $\frac{1}{r^3}$ singularity of $K^\text{MFIE}(r)$
at the origin might appear to preclude application of the
TD reduction method for this kernel in the common-tetrahedron
case (for which, as noted in the previous section, we are
guaranteed only desingularization of kernels with 
singularities of $\frac{1}{r^2}$ or weaker). However, the 
vanishing of $P\sups{MFIE}$ at $\vb x=\vb x^\prime$
affords extra leeway
(by ensuring that the sum over $n$ in (\ref{FinalIntegral})
begins at $n=1$),
allowing the pairing $\{P\sups{MEFIE}, K\sups{MEFIE}\}$
to be desingularized with no difficulty in the
common-tetrahedron and all other cases.

\subsection{Computation of VEFIE-SWG matrix elements: Comparison
            to Method of Bleszynski et al.}
\label{ComputationalResultsSection}
%####################################################################%
\begin{figure}
\begin{center}
\resizebox{0.5\textwidth}{!}{\includegraphics{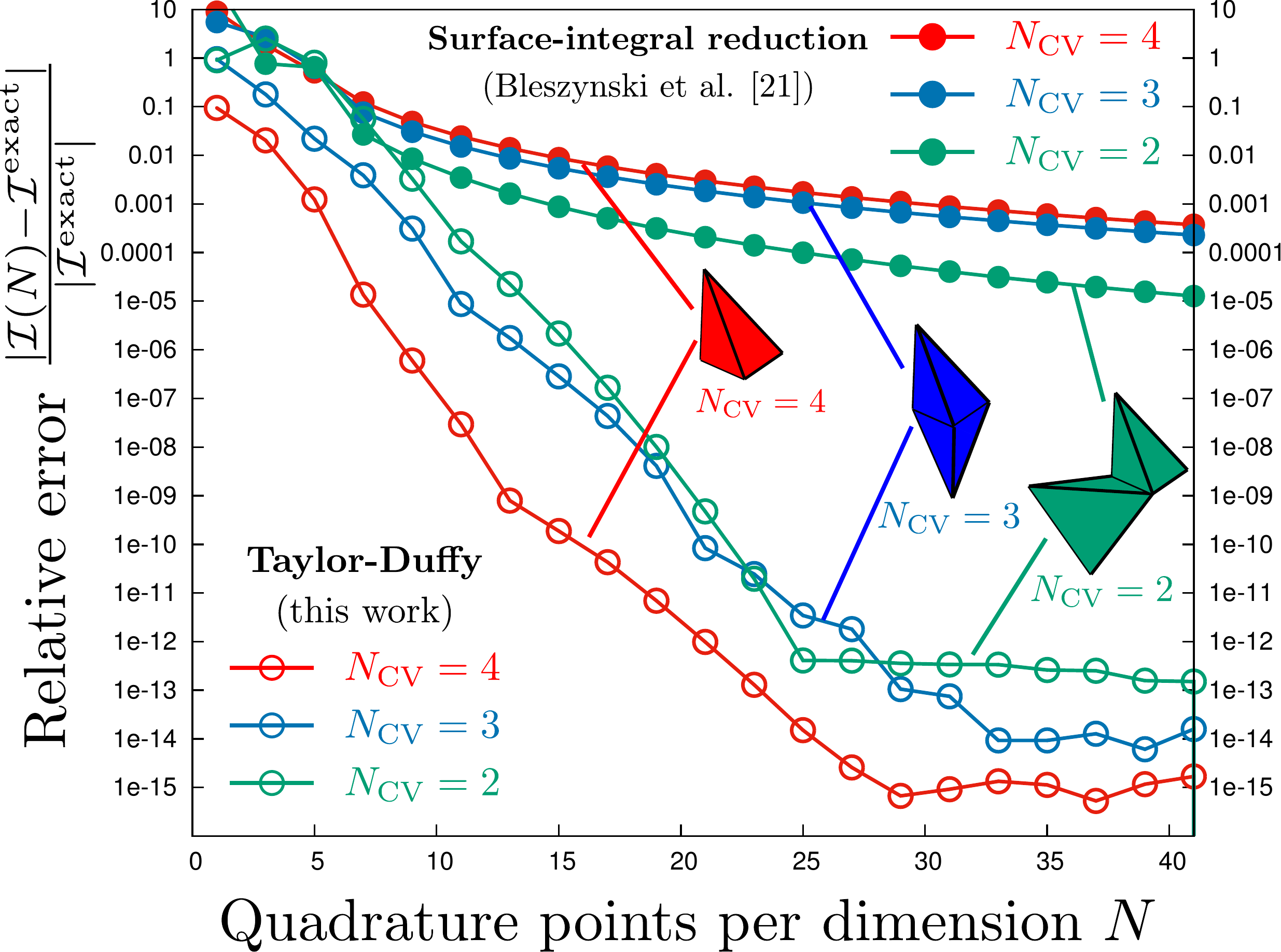}}
\caption{Comparison of convergence rates for the Taylor-Duffy (TD)
method proposed in this paper and for the surface-integral-reduction
(SIR) scheme of Bleszynski et al.~\cite{Bleszynski2013} as used
to evaluate VEFIE matrix elements [equation (\ref{OriginalIntegral})
with the $P$ and $K$ functions of equation (\ref{PKEFIE})]
at wavenumber $k=10$ in units where the tetrahedron edge lengths are of
order $\sim 1.$
For pairs of tetrahedra with $N\subt{CV}=\{4,3,2\}$ common
vertices (inset; see also Table \ref{TetVertexTable}),
I use the TD and SIR methods to transform the singular 6-dimensional 
integral (\ref{OriginalIntegral}) into a nonsingular $M$-dimensional 
integral (where $M=4$ for the SIR method and $M=6-N_\text{CV}=\{2,3,4\}$
for the TD method), then evaluate this integral numerically using nested
Clenshaw-Curtis quadrature with $N$ points per dimension 
(for a total of $N^M$ integrand samples) to
obtain an approximation $\mc I(N)$ to the original integral 
$(\ref{OriginalIntegral})$. 
Plotted is the relative error in this approximation vs. $N$.
The TD method converges exponentially with $N$, yielding
12 or more digits of accuracy for $N\approx 25$, and is
several orders of magnitude more accurate than the SIR 
method for all $N>10$.
For example, in the $N_\text{CV}=4$ case with $N=15$
the SIR method requires $15^4=50,625$ cubature
points to achieve 2-digit accuracy,
while the TD method requires $15^2=225$ points
to achieve 10-digit accuracy.
(As discussed in the text, the computational
cost per integrand sample is comparable for the two methods.)
}
\label{ConvergenceFigure}
\end{center}
\end{figure}
%####################################################################%

%####################################################################%
\begin{figure}
\begin{center}
\resizebox{0.5\textwidth}{!}{\includegraphics{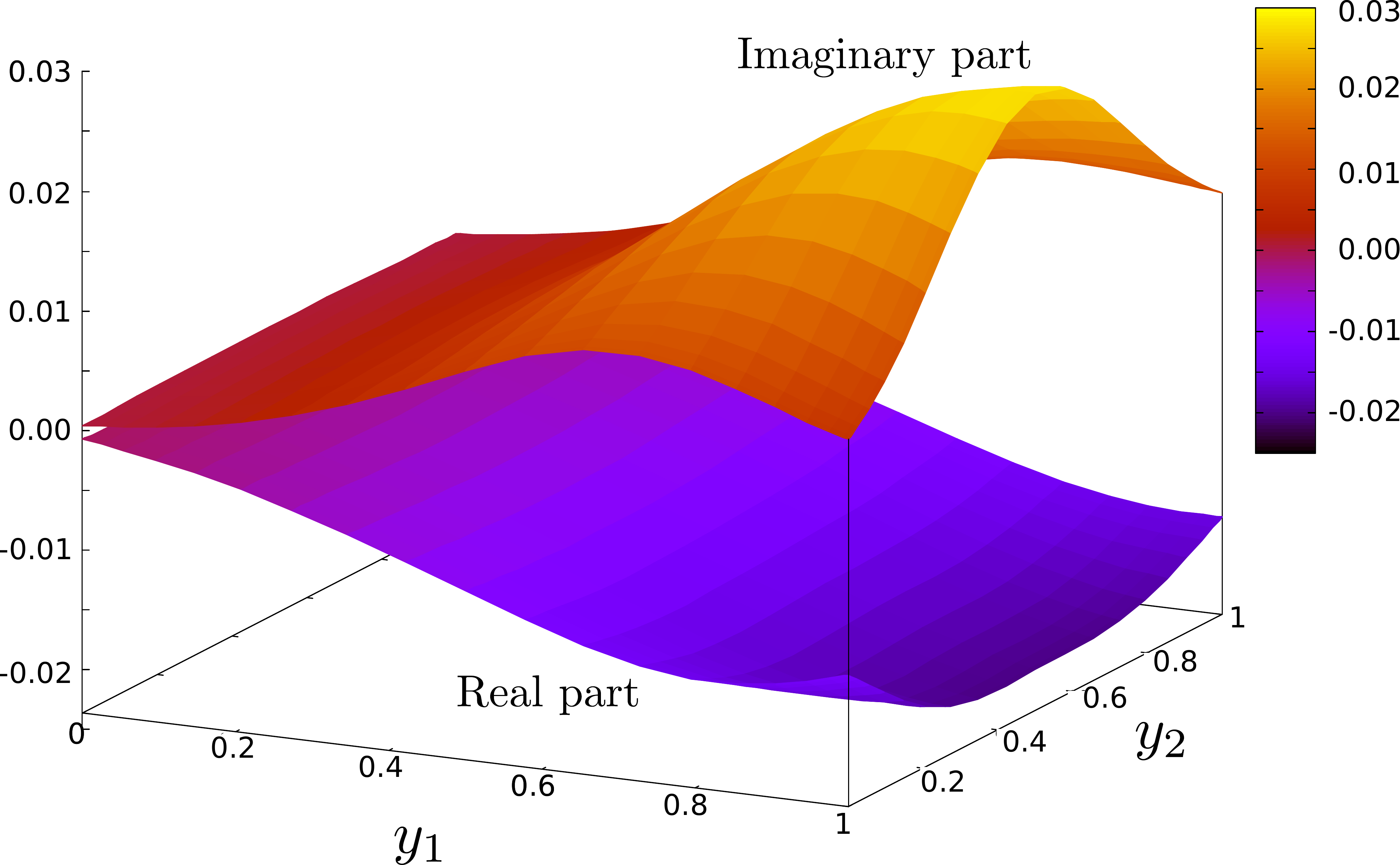}}
\caption{Integrand of the reduced integral (\ref{FinalIntegral})
for the common-tetrahedron ($N_\text{CV}=4$) case of Figure
\ref{ConvergenceFigure}. Whereas the integrand of the original
6-dimensional integral (\ref{OriginalIntegral}) has both integrable
singularities and sinusoidal variations at a wavelength 
$(\lambda=2\pi/10)$ shorter than the linear size $(L\sim 1)$ 
of the tetrahedron, the Taylor-Duffy reduction process
achieves a tremendous amount of smoothing; the gentle variation
of the reduced integrand with $\vb y$ explains why only
low-order cubature is required to evaluate the integral
to 15-digit accuracy (Figure \ref{ConvergenceFigure}).
}
\label{IntegrandFigure}
\end{center}
\end{figure}
%####################################################################%
As a concrete demonstration of the accuracy and efficiency of
the Taylor-Duffy (TD) method proposed here, I now use it to
compute the contributions of specific tetrahedron pairs to
VEFIE matrix elements between SWG basis functions---that is,
for fixed tetrahedra $\Tet,\Tet^\prime$ I evaluate equation
(\ref{OriginalIntegral}) with the $P$ and $K$ functions of
equation (\ref{PKEFIE}) at wavenumber $k=10$, corresponding
to a wavelength on the order of one-half the tetrahedron edge length.
For comparison, I also evaluate the same integrals using my
own implementation of the  
recently proposed surface-integral-reduction (SIR) method of
Bleszynski et al.~\cite{Bleszynski2013}.

I consider three pairs of tetrahedra with $N_\text{CV}=\{4,3,2\}$
common vertices (the common-\{tetrahedron, face, edge\} cases);
Table \ref{TetVertexTable} lists the vertices of these tetrahedra.
(Note that I have chosen tetrahedron $\mc T\subs{A}$ to be the
tetrahedron used by~\cite{Bleszynski2013} in SIR studies of
the $N_\text{CV}=4$ case, allowing direct comparison with results 
reported there.) The TD and SIR methods reduce the 6-dimensional 
singular integral (\ref{OriginalIntegral}) to a nonsingular 
$M$-dimensional 
integral (where $M=4$ for the SIR method and $M=6-N_\text{CV}=\{2,3,4\}$
for the TD method), which I evaluate numerically using nested
Clenshaw-Curtis (CC) quadrature~\cite{Trefethen2000} with $N$ 
points per dimension (total of $N^M$ integrand samples) to 
obtain an approximation
$\mc I(N)$ to the original integral (\ref{OriginalIntegral}).
(The SIR method involves integrals over triangles, which I
reparameterize as integrals over the unit square to allow
nested CC quadrature.)

Figure \ref{ConvergenceFigure} plots the relative error
$\mc E(N)\equiv |\mc I(N) - \mc I\sups{exact}|/|I\sups{exact}|$
versus $N$ for the TD and SIR methods.
(Reference values $\mc I\sups{exact}$ are the results
of TD calculations with $N=51$ and are tabulated in
Table \ref{IValuesTable}.)
My results for the SIR method in the $N_\text{CV}=4$ case
agree with the results of~\citeasnoun{Bleszynski2013},
which reported relative errors of $10^{-2}$ for $N=12$
and $10^{-3}$ for $N=40$. As the \textit{total} number
of cubature points used by the SIR method in these
cases are $N\subs{tot}=12^4\approx 2\cdot 10^4$ and 
$N\subs{tot}=40^4\approx 3\cdot 10^6$, the tenfold
reduction in error requires more than a hundredfold
increase in cost; the SIR method for this case appears
to be converging at the extremely slow algebraic rate
$\mc E\sim 1/\sqrt{N\subs{tot}}$ and is effectively
incapable of revealing more than a few correct digits
of $\mc I$ in practice.

In contrast, the TD method achieves \textit{exponential}
convergence in all cases, obtaining 12 or more correct
digits with as few as 25 cubature points per dimension.
For example, in the common-tetrahedron case with
$N=15$, the TD method requires a total of $N^2=225$
integrand samples to achieve 10-digit accuracy,
while the SIR method requires $15^4=50,625$ samples
to achieve roughly 2-digit accuracy.

%####################################################################%
\begin{table}
\renewcommand{\arraystretch}{1.5}
$$\begin{array}{|c|c|c|c|}\hline
         & \Tet_A  & \Tet_B          & \Tet_C                \\\hline
 \vb V_2 & (1,0,0) & (1,0,0)         & (0,0,1)               \\\hline
 \vb V_3 & (0,1,0) & (0,1,0)         & (-0.04, -1.09, -0.05) \\\hline
 \vb V_4 & (0,0,1) & (0.3,0.4,-1.03) & (0.3, -0.4, -1.09)    \\\hline
\end{array}$$
\renewcommand{\arraystretch}{1.0}
\caption{Vertices of tetrahedra used for sample calculations.
All tetrahedra have one vertex at the origin, $\vb V_1=(0,0,0)$.
The common-\{tetrahedron, face, edge\} cases 
($N_\text{CV}=\{4,3,2\}$) correspond to tetrahedron 
pairs 
$\{(\Tet_A,\Tet_A), (\Tet_A,\Tet_B), (\Tet_A,\Tet_C)\}.$
}
\label{TetVertexTable}
\end{table}
%####################################################################%

%####################################################################%
\begin{table}
\renewcommand{\arraystretch}{1.5}
$$\begin{array}{c|c|c|c}
 (\Tet, \Tet^\prime) & N_\text{CV} & 
 \text{Re }\mathcal{I}\sups{exact} &
 \text{Im }\mathcal{I}\sups{exact} 
\\\hline
 (\Tet_A, \Tet_A) & 4 & 
 \texttt{-7.8624620487335e-04} &
 \texttt{+8.5795441769385e-04}
\\ 
 (\Tet_A, \Tet_B) & 3 & 
 \texttt{+4.2568610165422e-05} &
 \texttt{+3.2199164645680e-05}
\\
 (\Tet_A, \Tet_C) & 2 & 
 \texttt{-3.0105189689052e-05} &
 \texttt{-7.1022045556570e-07}
\end{array}$$
\renewcommand{\arraystretch}{1.0}
\caption{Reference values of integral (\ref{OriginalIntegral})
with the $P$ and $K$ functions of equations (\ref{PKEFIE})
($k=10$), obtained by nested Clenshaw-Curtis 
quadrature of the reduced Taylor-Duffy integral (\ref{FinalIntegral})
with $N=51$ quadrature points per dimension.}
\label{IValuesTable}
\end{table}

\section{Conclusions}
\label{ConclusionsSection}

In this paper I extended the Taylor-Duffy approach to
singular Galerkin integrals, which had previously
been applied to integrals over triangle-product
domains~\cite{Taylor2003,Reid2015A} to the more 
challenging case of tetrahedron-product domains.
As I demonstrated, this yields an algorithm
for computing VIE matrix elements with accuracy and
efficiency exceeding those of existing methods by 
several orders of magnitude. I am hopeful that this 
new technique and its free-software
implementation~\cite{buff-em} will prove useful
for accelerating VIE solvers for electromagnetic
scattering and other physical applications.
Meanwhile, the successful extension from triangles
to tetrahedra testifies to the broad generality
of the basic Taylor-Duffy strategy---and suggests
that the full extent of its utility remains far from
fully explored.
%####################################################################%
\appendices
\newpage
\section{Tables of Subdomain-Dependent Quantities} 
\label{Appendix}

Tables \ref{TetTetBarTable1}-\ref{DuffyTransformTableNCV2}
provide explicit definitions of various subregion-dependent
quantities referenced in Section \ref{DerivationSection}.

\begin{table}
 $$ \begin{array}{|c|c|c|c|c|c|c|} \hline
    &&&&&&\\[-3pt]
%--------------------------------------------------------------------%
    d 
    & u_{1d}^{\text{\tiny{\sc min}}}
    & u_{1d}^{\text{\tiny{\sc max}}}
    & u_{2d}^{\text{\tiny{\sc min}}}
    & u_{2d}^{\text{\tiny{\sc max}}}
    & u_{3d}^{\text{\tiny{\sc min}}}
    & u_{3d}^{\text{\tiny{\sc max}}}
\\[4pt] \hline
%--------------------------------------------------------------------%
    1 & 0 & 1 & u_1 & 1 & u_2 & 1
\\\hline
%--------------------------------------------------------------------%
    2 & 0 & 1 & u_1 & 1 & 0 & u_2
\\\hline
%--------------------------------------------------------------------%
    3 & 0 & 1 & u_1 & 1 & u_2-1 & 0 
\\\hline
%--------------------------------------------------------------------%
    4 & 0 & 1 & 0 & u_1 & u_2 & 1-u_1+u_2
\\\hline
%--------------------------------------------------------------------%
    5 & 0 & 1 & 0 & u_1 & 0 & u_2 
\\\hline
%--------------------------------------------------------------------%
    6 & 0 & 1 & 0 & u_1 & u_1-1 & 0 
\\\hline
    7 & 0 & 1 & u_1-1 & 0 & 0 & 1-u_1 + u_2
\\\hline
    8 & 0 & 1 & u_1-1 & 0 & u_2 & 0
\\\hline
    9 & 0 & 1 & u_1 -1 & 0 & u_1-1 & u_2 
\\\hline
   10 & -1 & 0 & 0 & u_1+1 & u_2 & 1+u_1
\\\hline
   11 & -1 & 0 & 0 & u_1+1 & 0 & u_2
\\\hline
   12 & -1 & 0 & 0 & u_1+1 & u_2 - u_1 - 1 & 0
\\\hline
   13 & -1 & 0 & u_1 & 0 & 0 & u_1+1
\\\hline
   14 & -1 & 0 & u_1 & 0 & u_2 & 0
\\\hline
   15 & -1 & 0 & u_1 & 0 & u_2 - u_1 - 1 & u_2 
\\\hline
   16 & -1 & 0 & -1 & u_1 & 0 & 1+u_2
\\\hline
   17 & -1 & 0 & -1 & u_1 & u_2 &  0
\\\hline
   18 & -1 & 0 & -1 & u_1 & -1 & u_2
\\\hline
    \end{array}
%--------------------------------------------------------------------%
 $$
 \caption{ Limits of integration for the tetrahedral domains 
           $\tet^{\vb u}_d$ in (\ref{OriginalIntegral3}), defined by
           $\displaystyle{ 
             \int_{\tet_d} \, d\vb u
             \equiv \int_{u_{1d}\supt{min}}^{u_{1d}\supt{max}} \, du_1
                    \int_{u_{2d}\supt{min}}^{u_{2d}\supt{max}} \, du_2
                    \int_{u_{3d}\supt{min}}^{u_{3d}\supt{max}} \, du_3.
                         } 
           $
         }
 \label{TetTetBarTable1}
\end{table}

%%%%%%%%%%%%%%%%%%%%%%%%%%%%%%%%%%%%%%%%%%%%%%%%%%%%%%%%%%%%%%%%%%%%%%
%%%%%%%%%%%%%%%%%%%%%%%%%%%%%%%%%%%%%%%%%%%%%%%%%%%%%%%%%%%%%%%%%%%%%%
%%%%%%%%%%%%%%%%%%%%%%%%%%%%%%%%%%%%%%%%%%%%%%%%%%%%%%%%%%%%%%%%%%%%%%
\begin{table}
 $$ \begin{array}{|c|c|c|c|c|c|c|} \hline 
%--------------------------------------------------------------------%
    d 
    & L_{3d}
    & U_{3d}
    & L_{2d}
    & U_{2d}
    & L_{1d}
    & U_{1d}
\\ \hline
1 & 0     & -u_3 & u_3\!-\!u_2             & - u_2 & u_2\!-\!u_1 & - u_1 \\\hline
2 & 0     & -u_2 & 0                      & - u_2 & u_2\!-\!u_1 & - u_1 \\\hline 
3 & -u_3  & -u_2 & 0                      & - u_2 & u_2\!-\!u_1 & - u_1 \\\hline 
4 & 0     & u_2\!-\!u_1\!-\!u_3 & u_3\!-\!u_2 & - u_1 & 0          & - u_1 \\\hline
5 & 0     & -u_1 & 0                      & - u_1 & 0          & - u_1 \\\hline
6 & -u_3  & -u_1 & 0                      & - u_1 & 0          & - u_1 \\\hline
7 & 0     & u_2\!-\!u_1\!-\!u_3 & u_3\!-\!u_2 & - u_1 & 0          & - u_1 \\\hline
8 & -u_3  & u_2\!-\!u_1\!-\! u_3 & u_3\!-\!u_2 & - u_1 & 0          & - u_1 \\\hline
9 & -u_3  & -u_1 & 0                      & - u_1 & 0          & - u_1 \\\hline
10 & 0    & u_1 - u_3 & u_3\!-\!u_2        & u_1 - u_2 & u_2\!-\!u_1 & 0 \\\hline
11 & 0    & u_1 - u_2 & 0                 & u_1 - u_2 & u_2\!-\!u_1 & 0 \\\hline
12 & -u_3 & u_1 - u_2 & 0                 & u_1 - u_2 & u_2\!-\!u_1 & 0 \\\hline
13 & 0    & u_1 - u_3 & u_3\!-\!u_2        & u_1 - u_2 & u_2\!-\!u_1 & 0 \\\hline
14 & -u_3 & u_1 - u_3 & u_3\!-\!u_2        & u_1 - u_2 & u_2\!-\!u_1 & 0 \\\hline
15 & -u_3 & u_1 - u_2 & 0                 & u_1 - u_2 & u_2\!-\!u_1 & 0 \\\hline
16 & 0    & u_2 - u_3 & u_3\!-\!u_2        & 0 & 0 & 0 \\\hline
17 & -u_3 & u_2 - u_3 & u_3\!-\!u_2        & 0 & 0 & 0 \\\hline
18 & -u_3 & 0         & 0                 & 0 & 0 & 0 \\\hline
%--------------------------------------------------------------------%
\end{array}
 $$
 \caption{ Limits of integration for the tetrahedral domains
           $\tet^{\vbXi}_d$ in (\ref{OriginalIntegral3}), defined by
           $\displaystyle{ 
             \int_{\tet_d} \, d\vbXi
             \equiv \int_{L_{3d}}^{1+U_{3d}}         \, d\xi_3
                    \int_{\xi_3 + L_{2d}}^{1+U_{2d}} \, d\xi_2
                    \int_{\xi_2 + L_{1d}}^{1+U_{1d}} \, d\xi_1.
                         } 
           $
         }
 \label{TetTetBarTable2}
\end{table}
%%%%%%%%%%%%%%%%%%%%%%%%%%%%%%%%%%%%%%%%%%%%%%%%%%%%%%%%%%%%%%%%%%%%

\renewcommand{\arraystretch}{2}
\begin{table}
 $$ 
%--------------------------------------------------------------------%
%- left table, D=1..9
%--------------------------------------------------------------------%
\begin{array}{|c|c|c|c|c|}\hline
%--------------------------------------------------------------------%
    d 
    & \displaystyle{
        \mc J_d(\vb y)
                   }
    & \displaystyle{
       u_1(w,\vb y)
                   }
    & \displaystyle{
       u_2(w,\vb y)
                   }
    & \displaystyle{
       u_3(w,\vb y)
                   }
\\\hline
%--------------------------------------------------------------------%
    1 
    &   y_1   
    &  w y_1 y_2 
    &  w y_1     
    &  w         
\\\hline
%--------------------------------------------------------------------%
    2 
    &  1      
    &  w y_1     
    &  w         
    &  w y_2     
\\\hline
%--------------------------------------------------------------------%
    3 
    &   y_1   
    &  w y_1 y_2 
    &  w y_1     
    &  -w(1-y_1) 
\\\hline
%--------------------------------------------------------------------%
    4 
    &   y_1   
    &  w y_1     
    &  w y_1 y_2 
    &  w (1-y_1+y_1 y_2) 
\\\hline
%--------------------------------------------------------------------%
    5 
    &   y_1   
    &  w          
    &  w y_1     
    &  w y_1 y_2 
\\\hline
%--------------------------------------------------------------------%
    6 
    &   y_1   
    &  w y_1     
    &  w y_1 y_2 
    & -w(1-y_1)
\\\hline
%--------------------------------------------------------------------%
    7 
    &   y_1   
    &  wy_1y_2 
    &  -wy_1(1-y_2)
    & w(1-y_1) 
\\\hline
%--------------------------------------------------------------------%
    8 
    &   y_1   
    &  w(1-y_1)
    &  -wy_1
    &  -wy_1 y_2
\\\hline
%--------------------------------------------------------------------%
    9 
    &   y_1   
    &  wy_1(1-y_2) 
    &  -wy_1 y_2
    &  -w(1-y_1+y_1y_2)
\\\hline
% uncomment the next 23 lines to make the table double-wide and half-height
%\end{array}
%%--------------------------------------------------------------------%
%%- end of left table ------------------------------------------------%
%%--------------------------------------------------------------------%
%\qquad
%%--------------------------------------------------------------------%
%%- right table, d=10...18--------------------------------------------%
%%--------------------------------------------------------------------%
%\begin{array}{|c|c|c|c|c|}\hline
%%--------------------------------------------------------------------%
%    d 
%    & \displaystyle{
%        \mc J_d(w,\vb y)
%                   }
%    & \displaystyle{
%       u_1(w,\vb y)
%                   }
%    & \displaystyle{
%       u_2(w,\vb y)
%                   }
%    & \displaystyle{
%       u_3(w,\vb y)
%                   }
%\\\hline
%--------------------------------------------------------------------%
    10
    &   y_1   
    & -wy_1 y_2 
    &  wy_1(1-y_2)
    &  w(1-y_1 y_2)
\\\hline
%--------------------------------------------------------------------%
    11
    &    y_1   
    &  -wy_1 y_2
    &  w(1-y_1y_2)
    & wy_1(1-y_2)
\\\hline
%--------------------------------------------------------------------%
    12
    &    y_1     
    &  -wy_1 y_2
    &  wy_1 (1-y_2)
    & -w(1-y_1)
\\\hline
%--------------------------------------------------------------------%
    13
    &   y_1   
    & -wy_1
    & -wy_1 y_2
    & w(1-y_1)
\\\hline
%--------------------------------------------------------------------%
    14 
    &   y_1    
    & -w
    & -wy_1
    & -wy_1 y_2
\\\hline
%--------------------------------------------------------------------%
    15
    &   y_1   
    & -wy_1 
    & -wy_1 y_2
    & -w(1-y_1+y_1 y_2)
\\\hline
%--------------------------------------------------------------------%
    16
    &   y_1   
    & -wy_1 y_2 
    & -wy_1
    & w(1-y_1)
\\\hline
%--------------------------------------------------------------------%
    17
    & 1        
    & -wy_1 
    & -w 
    & -w(1-y_2)
\\\hline
%--------------------------------------------------------------------%
    18
    &   y_1   
    & -wy_1 y_2
    & -wy_1
    & -w
\\\hline
%--------------------------------------------------------------------%
\end{array}$$
 \caption{Duffy transformations for the case $N\subt{CV}=4$.} 
 \label{DuffyTableNCV4}
\end{table}
\renewcommand{\arraystretch}{1}

%%%%%%%%%%%%%%%%%%%%%%%%%%%%%%%%%%%%%%%%%%%%%%%%%%
%%%%%%%%%%%%%%%%%%%%%%%%%%%%%%%%%%%%%%%%%%%%%%%%%%
%%%%%%%%%%%%%%%%%%%%%%%%%%%%%%%%%%%%%%%%%%%%%%%%%%

\renewcommand{\arraystretch}{2}
\begin{table}
 $$ 
%--------------------------------------------------------------------%
\begin{array}{|c|c|c|c|c|c|}\hline
%--------------------------------------------------------------------%
    d 
    & \displaystyle{
        \mc J_d(\vb y)
                   }
    & \displaystyle{
       u_1(w,\vb y)
                   }
    & \displaystyle{
       u_2(w,\vb y)
                   }
    & \displaystyle{
       u_3(w,\vb y)
                   }
    & \displaystyle{
      \xi_3(w,\vb y)
                   }
\\\hline
%--------------------------------------------------------------------%
    1 
    &   y_1^2 y_2   
    &  w y_1 y_2 y_3   
    &  w y_1 y_2       
    &  w y_1            
    &  w \overline{y}_1 
\\\hline
%--------------------------------------------------------------------%
    2 
    &   y_1^2    
    &  w y_1 y_2 
    &  w y_1      
    &  w y_1 y_3  
    &  w \overline{y}_1  
\\\hline
%--------------------------------------------------------------------%
    3 
    &   y_1^2 y_2 
    &  w y_1 y_2 y_3 
    &  w y_1 y_2     
    &  -w y_1\overline{y}_2       
    &  w \zeta\subt{A}
\\\hline
%--------------------------------------------------------------------%
    4 
    &   y_1^2 y_2         
    &  w y_1 y_2             
    &  w y_1 y_2 y_3          
    &  w y_1 \zeta\subt{D}
    &  w \overline{y}_1              
\\\hline
%--------------------------------------------------------------------%
    5 
    &   y_1^2 y_2         
    &  w y_1
    &  w y_1 y_2
    &  w y_1 y_2 y_3
    &  w \overline{y}_1              
\\\hline
%--------------------------------------------------------------------%
    6 
    &   y_1^2 y_2         
    &  w y_1 y_2             
    &  w y_1 y_2 y_3         
    &  -w y_1 \overline{y}_2
    &  w \zeta\subt{A}
\\\hline
%--------------------------------------------------------------------%
    7 
    &   y_1^2 y_2         
    &  wy_1 y_2 y_3
    &  -wy_1 y_2\overline{y}_3
    &  wy_1\overline{y}_2
    &  w\overline{y}_1
\\\hline
%--------------------------------------------------------------------%
    8 
    &   y_1^2 y_2         
    &  wy_1\overline{y}_2
    &  -wy_1 y_2
    &  -wy_1 y_2 y_3
    &  w\zeta\subt{E}
\\\hline
%--------------------------------------------------------------------%
    9 
    &   y_1^2 y_2         
    &  wy_1 y_2\overline{y}_3
    &  -wy_1 y_2 y_3 
    &  -wy_1\zeta\subt{D}
    &  w\zeta\subt{F}
\\\hline
%--------------------------------------------------------------------%
   10 
    &   y_1^2 y_2         
    & -wy_1 y_2 y_3
    & wy_1 y_2\overline{y}_3
    & wy_1\zeta\subt{C}
    & w\overline{y}_1
\\\hline
%--------------------------------------------------------------------%
   11 
    &   y_1^2 y_2         
    & -wy_1 y_2 y_3
    & wy_1\zeta\subt{C}
    & wy_1 y_2\overline{y}_3
    & w\overline{y}_1
\\\hline
%--------------------------------------------------------------------%
   12 
    &   y_1^2 y_2         
    & -wy_1 y_2 y_3
    & wy_1 y_2\overline{y}_3
    & -wy_1\overline{y}_2
    & w\zeta\subt{A}
\\\hline
%--------------------------------------------------------------------%
   13
    &   y_1^2 y_2
    & -wy_1 y_2
    & -wy_1 y_2 y_3
    &  wy_1\overline{y}_2
    &  w\overline{y}_1
\\\hline
%--------------------------------------------------------------------%
   14 
    &   y_1^2 y_2         
    & -wy_1
    & -wy_1 y_2
    & -wy_1 y_2 y_3
    &  w\zeta\subt{D}
\\\hline
%--------------------------------------------------------------------%
   15 
    &   y_1^2 y_2         
    & -wy_1 y_2
    & -wy_1 y_2 y_3
    & -wy_1\zeta\subt{D}
    & w\zeta\subt{F}
\\\hline
%--------------------------------------------------------------------%
   16 
    &   y_1^2 y_2         
    & -wy_1 y_2 y_3
    & -wy_1 y_2
    & wy_1\overline{y}_2
    & w\overline{y}_1
\\\hline
%--------------------------------------------------------------------%
   17 
    &   y_1^2
    & -wy_1 y_2
    & -wy_1
    & -wy_1 \overline{y}_3
    & w\zeta\subt{B}
\\\hline
%--------------------------------------------------------------------%
   18 
    &   y_1^2 y_2         
    & -wy_1 y_2 y_3
    & -wy_1 y_2 
    & -wy_1 
    & w
\\\hline
%--------------------------------------------------------------------%
\end{array}$$
 \caption{Duffy transformations for the case $N\subt{CV}=3$.
          Shorthand: $\overline{y_i}\equiv 1-y_i$,\,
                     $\zeta\subt{A}\equiv 1-y_1 y_2$,\,
                     $\zeta\subt{B}\equiv 1-y_1 y_3$,\,
                     $\zeta\subt{C}\equiv 1-y_2 y_3$,\,
                     $\zeta\subt{D}\equiv 1-y_2+y_2 y_3$,\,
                     $\zeta\subt{E}\equiv 1-y_1+y_1y_2y_3$,\,
                     $\zeta\subt{F}\equiv 1-y_1y_2+y_1y_2y_3$.}
 \label{DuffyTableNCV3}
\end{table}
\renewcommand{\arraystretch}{1}

%%%%%%%%%%%%%%%%%%%%%%%%%%%%%%%%%%%%%%%%%%%%%%%%%%
%%%%%%%%%%%%%%%%%%%%%%%%%%%%%%%%%%%%%%%%%%%%%%%%%%
%%%%%%%%%%%%%%%%%%%%%%%%%%%%%%%%%%%%%%%%%%%%%%%%%%
\renewcommand{\arraystretch}{2}
\begin{table}
\hspace{-0.5in}
{\scriptsize
 $$ 
%--------------------------------------------------------------------%
\begin{array}{|c|c|c|c|c|c|c|}\hline
%--------------------------------------------------------------------%
    d 
    & \displaystyle{
        \mc J_d
                   }
    & \displaystyle{
       \frac{\vphantom{\int} u_1}{w}
                   }
    & \displaystyle{
       \frac{u_2}{w}
                   }
    & \displaystyle{
       \frac{u_3}{w}
                   }
    & \displaystyle{
      \frac{\xi_3}{w}
                   }
    & \displaystyle{
      \frac{\xi_2}{w}
                   }
\\[0.05in]\hline
%--------------------------------------------------------------------%
    1 
    &   y_1^3 y_2^2 y_3   
    &   y_1 y_2 y_3 y_4     
    &   y_1 y_2 y_3         
    &   y_1 y_2             
    &   y_1\overline{y}_2          
    &   \Upsilon\subt{E}
\\\hline
%--------------------------------------------------------------------%
    2 
    &   y_1^3 y_2^2       
    &   y_1 y_2 y_4         
    &   y_1 y_2             
    &   y_1 y_2 y_3         
    &   y_1\overline{y}_2          
    &   \Upsilon\subt{A}
\\\hline
%--------------------------------------------------------------------%
    3 
    &   y_1^3 y_2^2 y_3   
    &   y_1 y_2 y_3 y_4     
    &   y_1 y_2 y_3         
    &  - y_1 y_2 \overline{y}_3            
    &  y_1 \Upsilon\subt{B}
    &   \Upsilon\subt{E}
\\\hline
%--------------------------------------------------------------------%
    4 
    &   y_1^3 y_2^2 y_3   
    &   y_1 y_2 y_3         
    &   y_1 y_2 y_3 y_4     
    &   y_1 y_2 \Upsilon\subt{G}
    &   y_1 \overline{y}_2       
    &   \Upsilon\subt{E}
\\\hline
%--------------------------------------------------------------------%
    5 
    &   y_1^3 y_2^2 y_3   
    &   y_1 y_2
    &   y_1 y_2 y_3
    &   y_1 y_2 y_3 y_4
    &   y_1 \overline{y}_2       
    &   \Upsilon\subt{A}
\\\hline
%--------------------------------------------------------------------%
    6 
    &   y_1^3 y_2^2 y_3   
    &   y_1 y_2 y_3         
    &   y_1 y_2 y_3 y_4     
    &  - y_1 y_2\overline{y_3}
    &   y_1 \Upsilon\subt{B}
    &  \Upsilon\subt{E}
\\\hline
%--------------------------------------------------------------------%
    7 
    &   y_1^3 y_2^2 y_3   
    & y_1 y_2 y_3 y_4
    & - y_1 y_2 y_3\overline{y}_4
    & y_1 y_2\overline{y}_3
    & y_1\overline{y}_2
    & \Upsilon\subt{F}
\\\hline
%--------------------------------------------------------------------%
    8 
    &   y_1^3 y_2^2 y_3   
    & y_1y_2\overline{y}_3
    & -y_1y_2y_3
    & -y_1y_2y_3y_4
    & y_1\Upsilon\subt{H}
    & \Upsilon\subt{I}
\\\hline
%--------------------------------------------------------------------%
    9 
    &   y_1^3 y_2^2 y_3   
    & y_1y_2y_3\overline{y}_4
    &  -y_1 y_2 y_3 y_4
    & -y_1y_2\Upsilon\subt{G}
    & y_1\Upsilon\subt{J}
    & \Upsilon\subt{K}
\\\hline
%--------------------------------------------------------------------%
   10 
    &   y_1^3 y_2^2 y_3   
    & -y_1 y_2 y_3 y_4
    & y_1 y_2 y_3\overline{y}_4
    & y_1 y_2\Upsilon\subt{D}
    & y_1\overline{y}_2
    & \Upsilon\subt{E}
\\\hline
%--------------------------------------------------------------------%
   11 
    &   y_1^3 y_2^2 y_3   
    & -y_1 y_2 y_3 y_4
    & y_1 y_2\Upsilon\subt{D}
    & y_1y_2y_3\overline{y}_4
    & y_1\overline{y}_2
    & \Upsilon\subt{A}
\\\hline
%--------------------------------------------------------------------%
   12 
    &   y_1^3 y_2^2 y_3   
    & -y_1 y_2 y_3 y_4
    & y_1y_2y_3\overline{y}_4
    & -y_1y_2\overline{y}_3
    & y_1\Upsilon\subt{B}
    & \Upsilon\subt{E}
\\\hline
%--------------------------------------------------------------------%
   13 
    &   y_1^3 y_2^2 y_3   
    &  -y_1y_2y_3
    &  -y_1y_2y_3y_4
    &   y_1y_2\overline{y}_3
    &   y_1\overline{y}_2
    &   \Upsilon\subt{K}
\\\hline
%--------------------------------------------------------------------%
   14 
    &   y_1^3 y_2^2 y_3   
    & -y_1 y_2
    & -y_1 y_2 y_3
    & -y_1 y_2 y_3 y_4
    & y_1\Upsilon\subt{H}
    & \Upsilon\subt{I}
\\\hline
%--------------------------------------------------------------------%
   15 
    &   y_1^3 y_2^2 y_3   
    & -y_1 y_2 y_3
    & -y_1 y_2 y_3 y_4
    & -y_1 y_2\Upsilon\subt{G}
    & y_1\Upsilon\subt{J}
    & \Upsilon\subt{K}
\\\hline
%--------------------------------------------------------------------%
   16 
    &   y_1^3 y_2^2 y_3   
    &  -y_1 y_2 y_3 y_4
    & -y_1 y_2 y_3
    & y_1 y_2\overline{y}_3
    & y_1 \overline{y}_2
    & 1
\\\hline
%--------------------------------------------------------------------%
   17 
    &   y_1^3 y_2^2 y_3   
    &  -y_1 y_2 y_3
    & -y_1 y_2
    & -y_1y_2\overline{y}_4
    & y_1\Upsilon\subt{C}
    & 1
\\\hline
%--------------------------------------------------------------------%
   18 
    &   y_1^3 y_2^2 y_3   
    & -y_1 y_2 y_3 y_4
    & -y_1 y_2 y_3
    & -y_1 y_2
    & y_1
    & 1
\\\hline
%--------------------------------------------------------------------%
\end{array}$$
} % { \scriptsize
 \caption{Duffy transformations for the case $N\subt{CV}=2$.
          (Note that a factor of $w$ has been extracted from
          each $\vb u, \vbXi$ component to save space).
          Shorthand: $\overline{y}_i \equiv 1-y_i$,\, 
                     $\Upsilon\subt{A}\equiv 1-y_1y_2$,\,
                     $\Upsilon\subt{B}\equiv 1-y_2y_3$,\,
                     $\Upsilon\subt{C}\equiv 1-y_2y_4$,\,
                     $\Upsilon\subt{D}\equiv 1-y_3y_4$,\,
                     $\Upsilon\subt{E}\equiv 1-y_1y_2y_3$,\,
                     $\Upsilon\subt{F}\equiv 1-y_1y_2y_3y_4$,\,
                     $\Upsilon\subt{G}\equiv 1-y_3+y_3y_4$,\,
                     $\Upsilon\subt{H}\equiv 1-y_2+y_2y_3y_4$,\,
                     $\Upsilon\subt{I}\equiv 1-y_1y_2+y_1y_2y_3$,\,
                     $\Upsilon\subt{J}\equiv 1-y_2y_3+y_2y_3y_4$,\,
                     $\Upsilon\subt{K}\equiv 1-y_1y_2y_3+y_1y_2y_3y_4$.} \label{DuffyTableNCV2}
\label{DuffyTransformTableNCV2}
\end{table}
\renewcommand{\arraystretch}{1}

%%%%%%%%%%%%%%%%%%%%%%%%%%%%%%%%%%%%%%%%%%%%%%%%%%%%%%%%%%%%%%%%%%%%%%%%%%%
%%%%%%%%%%%%%%%%%%%%%%%%%%%%%%%%%%%%%%%%%%%%%%%%%%%%%%%%%%%%%%%%%%%%%%%%%%%
%%%%%%%%%%%%%%%%%%%%%%%%%%%%%%%%%%%%%%%%%%%%%%%%%%%%%%%%%%%%%%%%%%%%%%%%%%%
\bibliographystyle{IEEEtran}

\begin{thebibliography}{10}
\providecommand{\url}[1]{#1}
\csname url@samestyle\endcsname
\providecommand{\newblock}{\relax}
\providecommand{\bibinfo}[2]{#2}
\providecommand{\BIBentrySTDinterwordspacing}{\spaceskip=0pt\relax}
\providecommand{\BIBentryALTinterwordstretchfactor}{4}
\providecommand{\BIBentryALTinterwordspacing}{\spaceskip=\fontdimen2\font plus
\BIBentryALTinterwordstretchfactor\fontdimen3\font minus
  \fontdimen4\font\relax}
\providecommand{\BIBforeignlanguage}[2]{{%
\expandafter\ifx\csname l@#1\endcsname\relax
\typeout{** WARNING: IEEEtran.bst: No hyphenation pattern has been}%
\typeout{** loaded for the language `#1'. Using the pattern for}%
\typeout{** the default language instead.}%
\else
\language=\csname l@#1\endcsname
\fi
#2}}
\providecommand{\BIBdecl}{\relax}
\BIBdecl

\bibitem{Harrington93}
R.~F. Harrington, \emph{Field Computation by Moment Methods}.\hskip 1em plus
  0.5em minus 0.4em\relax Wiley-IEEE Press, 1993.

\bibitem{Chew2009}
\BIBentryALTinterwordspacing
W.~Chew, M.~Tong, and B.~Hu, \emph{Integral Equation Methods for
  Electromagnetic and Elastic Waves}, ser. Synthesis Lectures on Computational
  Electromagnetics Series.\hskip 1em plus 0.5em minus 0.4em\relax Morgan \&
  Claypool Publishers, 2009. [Online]. Available:
  \url{http://books.google.com/books?id=PJN9meadzT8C}
\BIBentrySTDinterwordspacing

\bibitem{Volakis2012}
\BIBentryALTinterwordspacing
S.~K. Volakis, John.\hskip 1em plus 0.5em minus 0.4em\relax SciTech Publishing,
  2012. [Online]. Available:
  \url{http://app.knovel.com/hotlink/toc/id:kpIEME0011/integral-equation-methods/integral-equation-methods}
\BIBentrySTDinterwordspacing

\bibitem{Bleszynski2008}
E.~Bleszynski, M.~Bleszynski, and T.~Jaroszewicz, ``Fast volumetric
  integral-equation solver for acoustic wave propagation through inhomogeneous
  media,'' \emph{The Journal of the Acoustical Society of America}, vol. 124,
  no.~1, 2008.

\bibitem{Markkanen2014}
J.~Markkanen and P.~Ylä-Oijala, ``Discretization of electric current volume
  integral equation with piecewise linear basis functions,'' \emph{IEEE
  Transactions on Antennas and Propagation}, vol.~62, no.~9, pp. 4877--4880,
  Sept 2014.

\bibitem{Schaubert1984}
D.~Schaubert, D.~Wilton, and A.~Glisson, ``A tetrahedral modeling method for
  electromagnetic scattering by arbitrarily shaped inhomogeneous dielectric
  bodies,'' \emph{IEEE Transactions on Antennas and Propagation}, vol.~32,
  no.~1, pp. 77--85, Jan 1984.

\bibitem{Duffy1982}
M.~G. Duffy, ``Quadrature over a pyramid or cube of integrands with a
  singularity at a vertex,'' \emph{SIAM Journal on Numerical Analysis},
  vol.~19, no.~6, pp. 1260--1262, 1982.

\bibitem{Taylor2003}
D.~Taylor, ``Accurate and efficient numerical integration of weakly singular
  integrals in {G}alerkin {EFIE} solutions,'' \emph{Antennas and Propagation,
  IEEE Transactions on}, vol.~51, no.~7, pp. 1630--1637, 2003.

\bibitem{Reid2015A}
M.~Reid, J.~White, and S.~Johnson, ``Generalized Taylor-Duffy method for
  efficient evaluation of galerkin integrals in boundary-element method
  computations,'' \emph{Antennas and Propagation, IEEE Transactions on},
  vol.~63, no.~1, pp. 195--209, Jan 2015.

\bibitem{Bleszynski2013}
E.~Bleszynski, M.~Bleszynski, and T.~Jaroszewicz, ``Reduction of volume
  integrals to nonsingular surface integrals for matrix elements of tensor and
  vector green functions of maxwell equations,'' \emph{Antennas and
  Propagation, IEEE Transactions on}, vol.~61, no.~7, pp. 3642--3647, July
  2013.

\bibitem{Hasanovic2007}
M.~Hasanovic, C.~Mei, J.~R. Mautz, and E.~Arvas, ``Scattering from 3-d
  inhomogeneous chiral bodies of arbitrary shape by the method of moments,''
  \emph{IEEE Transactions on Antennas and Propagation}, vol.~55, no.~6, pp.
  1817--1825, June 2007.

\bibitem{Zhang2015}
L.~M. Zhang and X.~Q. Sheng, ``Solving volume electric current integral
  equation with full- and half-SWG functions,'' \emph{IEEE Antennas and
  Wireless Propagation Letters}, vol.~14, pp. 682--685, 2015.

\bibitem{Jackman2016}
\BIBentryALTinterwordspacing
K.~Jackman and C.~Fourie, ``Tetrahedral modeling method for inductance
  extraction of complex 3-d superconducting structures,'' \emph{IEEE
  Transactions on Applied Superconductivity}, vol.~26, no.~3, 2016, cited By 0.
  [Online]. Available:
  \url{https://www.scopus.com/inward/record.uri?eid=2-s2.0-84963795087&partnerID=40&md5=0c1d22abfe016803c1242966982aef81}
\BIBentrySTDinterwordspacing

\bibitem{Reid2016B}
M.~T.~H. Reid, ``Efficient Computation of Power, Force and Torque in
  Integral-Equation Solvers: Nonsingular Integrals and Moment Expansions'', to
  appear.

\bibitem{Botha2006}
\BIBentryALTinterwordspacing
M.~M. Botha, ``Solving the volume integral equations of electromagnetic
  scattering,'' \emph{Journal of Computational Physics}, vol. 218, no.~1, pp.
  141 -- 158, 2006. [Online]. Available:
  \url{http://www.sciencedirect.com/science/article/pii/S0021999106000763}
\BIBentrySTDinterwordspacing

\bibitem{Markkanen2012}
J.~Markkanen, C.~C. Lu, X.~Cao, and P.~Yla-Oijala, ``Analysis of volume
  integral equation formulations for scattering by high-contrast penetrable
  objects,'' \emph{IEEE Transactions on Antennas and Propagation}, vol.~60,
  no.~5, pp. 2367--2374, May 2012.

\bibitem{Polimeridis2014}
\BIBentryALTinterwordspacing
A.~Polimeridis, J.~Villena, L.~Daniel, and J.~White, ``Stable FFT-JVIE solvers
  for fast analysis of highly inhomogeneous dielectric objects,'' \emph{Journal
  of Computational Physics}, vol. 269, pp. 280 -- 296, 2014. [Online].
  Available:
  \url{http://www.sciencedirect.com/science/article/pii/S0021999114002071}
\BIBentrySTDinterwordspacing

\bibitem{Cools2003}
\BIBentryALTinterwordspacing
R.~Cools, ``An encyclopaedia of cubature formulas,'' \emph{Journal of
  Complexity}, vol.~19, no.~3, pp. 445 -- 453, 2003, oberwolfach Special Issue.
  [Online]. Available:
  \url{http://www.sciencedirect.com/science/article/pii/S0885064X03000116}
\BIBentrySTDinterwordspacing

\bibitem{RWG1982}
S.~Rao, D.~Wilton, and A.~Glisson, ``Electromagnetic scattering by surfaces of
  arbitrary shape,'' \emph{Antennas and Propagation, IEEE Transactions on},
  vol.~30, no.~3, pp. 409--418, May 1982.

\bibitem{Medgyesi1994}
\BIBentryALTinterwordspacing
L.~N. Medgyesi-Mitschang, J.~M. Putnam, and M.~B. Gedera, ``Generalized method
  of moments for three-dimensional penetrable scatterers,'' \emph{J. Opt. Soc.
  Am. A}, vol.~11, no.~4, pp. 1383--1398, Apr 1994. [Online]. Available:
  \url{http://josaa.osa.org/abstract.cfm?URI=josaa-11-4-1383}
\BIBentrySTDinterwordspacing

\bibitem{Hu2008}
L.~Hu, L.-W. Li, T.-S. Yeo, and R.~Vahldieck, ``An accurate and robust approach
  for evaluating vie impedance matrix elements using SWG basis functions,'' in
  \emph{Microwave Conference, 2008. APMC 2008. Asia-Pacific}, Dec 2008, pp.
  1--4.

\bibitem{Jarvenpaa2003}
\BIBentryALTinterwordspacing
S.~Järvenpää, M.~Taskinen, and P.~Ylä-Oijala, ``Singularity extraction
  technique for integral equation methods with higher order basis functions on
  plane triangles and tetrahedra,'' \emph{International Journal for Numerical
  Methods in Engineering}, vol.~58, no.~8, pp. 1149--1165, 2003. [Online].
  Available: \url{http://dx.doi.org/10.1002/nme.810}
\BIBentrySTDinterwordspacing

\bibitem{Wilton1984}
D.~Wilton, S.~Rao, A.~Glisson, D.~Schaubert, O.~Al-Bundak, and C.~Butler,
  ``Potential integrals for uniform and linear source distributions on
  polygonal and polyhedral domains,'' \emph{IEEE Transactions on Antennas and
  Propagation}, vol.~32, no.~3, pp. 276--281, Mar 1984.

\bibitem{DFN1985}
\BIBentryALTinterwordspacing
R.~Burns, B.~Dubrovin, A.~Fomenko, and S.~Novikov, \emph{Modern Geometry—
  Methods and Applications: Part II: The Geometry and Topology of Manifolds},
  ser. Graduate Texts in Mathematics.\hskip 1em plus 0.5em minus 0.4em\relax
  Springer New York, 1985. [Online]. Available:
  \url{https://books.google.com/books?id=tlzc7xXYKd8C}
\BIBentrySTDinterwordspacing

\bibitem{buff-em}
\texttt{https://github.com/HomerReid/buff-em}.

\bibitem{Helsinki2014}
J.~Markkanen and P.~Ylä-Oijala, ``Discretization of electric current volume
  integral equation with piecewise linear basis functions,'' \emph{IEEE
  Transactions on Antennas and Propagation}, vol.~62, no.~9, pp. 4877--4880,
  Sept 2014.

\bibitem{Helsinki2012}
J.~Markkanen, C.~C. Lu, X.~Cao, and P.~Yla-Oijala, ``Analysis of volume
  integral equation formulations for scattering by high-contrast penetrable
  objects,'' \emph{IEEE Transactions on Antennas and Propagation}, vol.~60,
  no.~5, pp. 2367--2374, May 2012.

\bibitem{Trefethen2000}
\BIBentryALTinterwordspacing
L.~Trefethen, \emph{Spectral Methods in MATLAB}, ser. Software, Environments,
  and Tools.\hskip 1em plus 0.5em minus 0.4em\relax Society for Industrial and
  Applied Mathematics, 2000. [Online]. Available:
  \url{https://books.google.cz/books?id=pB4xiZKZ4ecC}
\BIBentrySTDinterwordspacing

\end{thebibliography}


\end{document}